\documentclass [12pt] {amsart}
\usepackage{amsmath, amsfonts, latexsym, array, amssymb, }
\usepackage{psfrag}
\usepackage[all]{xypic}
\usepackage{graphicx}

\newtheorem*{thma}{Theorem A}
\newtheorem*{thmb}{Theorem B}

\newtheorem{thm}{Theorem}[section]
\newtheorem{cor}[thm]{Corollary}
\newtheorem{lem}[thm]{Lemma}
\newtheorem{prop}[thm]{Proposition}
\newtheorem{cla}[thm]{Claim}
\newtheorem{defn}[thm]{Definition}
\newtheorem{rem}[thm]{Remark}
\newtheorem{conj}[thm]{Conjecture}
\newtheorem{que}[thm]{Question}
\newtheorem{fac}[thm]{Fact}

\newenvironment{demo}[1][Proof]{{\smallskip
\noindent\bf #1:
}}{\hfill$\Box$\smallskip}

\thanks{This paper was partially supported
CNPq, Faperj, and PRONEX-Dynamical Systems (Brazil).
The authors thank BYU and PUC-Rio for the support and hospitality during their visits while preparing this paper.}

\begin{document}

\title[Symbolic extensions and transitivity] {Symbolic extensions for partially hyperbolic diffeomorphisms}
\author{Lorenzo J. D\'{i}az, Todd Fisher}
\address{Depto. Matem\'{a}tica, PUC-Rio, Marqu\^{e}s de S. Vicente 225,
22453-900 Rio de Janeiro RJ, Brazil}
\address{Department of Mathematics, Brigham Young University, Provo, UT 84602}
\email{lodiaz@mat.puc-rio.br, tfisher@math.byu.edu}
\thanks{}

\date{May 5, 2009}
\commby{}

\begin{abstract}
We show that every partially hyperbolic diffeomorphism with a 1-dimensional center bundle has a principal symbolic extension.  On the other hand, we show there are no symbolic extensions $C^1$-generically among diffeomorphisms containing non-hyperbolic robustly transitive sets 
with center indecomposable bundle of dimension at least 2.
 \end{abstract}

\maketitle

\section{Introduction}

Symbolic dynamical systems arose as a tool to code complicated dynamical 
systems by the use of more tractable systems (shift spaces).
The existence of Markov partitions for hyperbolic systems provides an effective coding of these systems using symbolic dynamics.  A natural question, that we address in the present work, is to what extent a non-hyperoblic system can be codified using symbolic systems.

A dynamical system $(X,f)$ has a {\it symbolic extension} if there exists a subshift $(Y,\sigma)$ and a continuous surjective map $\pi:Y\rightarrow X$ such that $\pi \circ\sigma=f\circ\pi$: the system $(Y,\sigma)$ is called an {\it extension} of $(X,f)$ and $(X,f)$ is called a {\it factor} of $(Y,\sigma)$.

\subsection{Existence of symbolic extensions}
Every dynamical system with a symbolic extension has finite entropy.
J. Auslander asked if the converse holds.  Namely, that if every system with finite entropy has a symbolic extension.   In fact,
the existence of symbolic extensions seems to be related to the notion of asymptotically $h$-ex\-pan\-si\-ve systems, as defined by Misiurecicz~\cite{Mis} (in Section~\ref{s.background} we review all relevant notions in the 
introduction). 
A nice form of a symbolic extension is a {\it principal} one, that is an extension given by a factor map which preserves entropy for every invariant measure.
Boyle, D.~Fiebig, and U.~Fiebig~\cite{BFF}   were able to show the following:

\begin{thm}[\cite{BFF}]\label{t.symbolicextensions}
If $(X,f)$ is asymptotically $h$-expansive, then $f$ has a principal symbolic extension.
\end{thm}

Using a result of Buzzi~\cite{Buz}, that states that any $C^{\infty}$ diffeomorphism is asymptoically $h$-expansive
one gets the following: 

\begin{cor} [\cite{BFF,Buz}]
Every $C^{\infty}$ map of a compact manifold has a principal symbolic extension.
\end{cor}

This result shows connections between the regularity of the maps and the existence of symbolic extensions.

From now on we let $\mathrm{Diff}^r(M)$ be the space of $C^r$ diffeomorphisms of a closed compact manifold $M$ endowed
with the usual uniform topology. When $r=1$ we simply write $\mathrm{Diff}(M)$. 

\begin{conj}[Downarowicz and Newhouse~\cite{DN}] If $f\in\mathrm{Diff}^r(M)$ and $r\geq 2$, then $f$ has a symbolic extension.
\end{conj}
For instance, in a recent work Downarowicz and Maass~\cite{DM} have shown the above conjecture for maps of a closed interval.

Besides the regularity of the system a second ingredient for the existence of symbolic extensions seems to be hyperbolic-like properties.
This is precisely the idea behind our first result.

\begin{thma}\label{t.existence}
Every partially hyperbolic diffeomorphism with a $1$-di\-men\-sio\-nal center bundle is asymptotically $h$-expansive and therefore has a principal symbolic extension.
\end{thma}

The above result supports the principle of Pugh and Shub that ``a little hyperbolicity goes a long way''~\cite{PS} and the suggestion made to us by Burns~\cite{Burns} that 1-dimensional center partially hyperbolic systems in many ways behave like the hyperbolic ones.

One important class of weakly hyperbolic dynamical systems are the robustly transitive ones.
Given an open set $U$ and a diffeomorphism $f$, we let $\Lambda_f(U)=\bigcap_{n\in \mathbb{Z}} f^n (\bar U)$; we say that the set
$\Lambda_f(U)$ is {\it robustly transitive\/} if $\Lambda_f(U)\subset U$ and for every $g$ $C^1$-close to $f$ the set $\Lambda_g(U)$ is 
{\it transitive\/} (has a point with a forward dense orbit). 
When $U$ is the whole manifold we say $f$ is {\it robustly transitive.\/}

The first examples of 
robustly transitive sets were the hyperbolic basic sets (including transitive Anosov diffeomorphisms).
Examples of robustly non-hyperbolic transitive diffeomorphisms were constructed by  Shub~\cite{Shu1}, Ma\~{n}\'{e}~\cite{man78},
and others \cite{BD,BV}. Finally, the robustly transitive sets were
introduced in \cite{DPU} as
a generalization of these examples. 
\begin{thm}[\cite{BDP, DPU}] \label{t.BDPU} Every robustly transitive set has a dominated splitting.
\end{thm}

In fact, there is a ``finest" dominated splitting: that is a dominated splitting $E_1\oplus\cdots \oplus E_{k(f)}$ of $T\Lambda_f(U)$ such that each subbundle $E_i$ of the splitting is indecomposibe (i.e.  there is no dominated splitting $G_i\oplus F_i$ of $E_i$), see~\cite[Theorem 4]{BDP}.
Additionally, some of the examples above are partially hyperbolic and have a 1-dimensional center bundle. We have the following corollary to Theorem~A.

\begin{cor}\label{c.existence}
Every robustly transitive set with a partially hyperbolic splitting and a 
center 1-dimensional direction has a symbolic extension.
\end{cor}

As a remark, we observe that this result shows that super-ex\-po\-nen\-tial growth of
the number of periodic points does not preclude the existence of symbolic extensions. Specifically,
the $C^1$ generic diffeomorphisms satisfying Corollary~\ref{c.existence}
have super-exponential growth of the number of periodic points \cite{BDF}.

\subsection{Non-eistence of symbolic extensions}
Boyle, D.~Fiebig, and U.~Fiebig~\cite{BFF} constructed a number of topological examples  that have no symbolic extensions. 
The first examples of diffeomorphisms with no symbolic extensions were constructed by Downarowicz and Newhouse~\cite{DN}.  They show that a generic area-preserving $C^1$ diffeomorphism of a surface is either Anosov or has no symbolic extension. 
This result relies on the existence of ``persistent homoclinic tangencies.''

 Asaoka~\cite{Asa}  gives a simple example of diffeomorphisms in a 3-di\-men\-sio\-nal disk with $C^1$ persistent 
homoclinic tangencies, and uses the result of Downarowicz and Newhouse, to show that for any smooth manifold $M$ with dim$(M) \geq 3$, there exists an open subset of $\mathrm{Diff}(M)$ in which generic diffeomorphisms have 
no symbolic extensions.  Finally, he extends this example to higher dimensions using normal hyperbolicity. 

We now describe another class of $C^1$ generic diffeomorphisms without symbolic extensions, and show how this is
related with non-hyperbolic properties.

Let $U$ be an open set in a closed manifold $M$.
Let $\mathcal{T}(U)$ be the set of diffeomorphisms, $f$, such that $\Lambda_f(U)$ is
robustly transitive.
Let $\mathcal{T}^{nh}(U)$ denote the set of  diffeomorphisms, $f$, such that
$\Lambda_f(U)$ is not hyperbolic.
We let $\mathcal{T}^{nh}_2(U)$ be the subset of $\mathcal{T}^{nh}(U)$ 
such that for each $f\in \mathcal{T}^{nh}_2(U)$ there is a neighborhood $\mathcal{V}$ of $f$ such that for each $g\in\mathcal{V}$ the set $\Lambda_g(U)$  has a non-hyperbolic center indecomposable bundle of dimension at least $2$.

\begin{thmb}\label{t.nonexistence}
There is a $C^1$ residual set $\mathcal{R}$ in $\mathcal{T}^{nh}_2(U)$ such that each diffeomorphism in $\mathcal{R}$ has no symbolic extension.
\end{thmb}

Relevant  systems satisfying this theorem are the DA-diffeomorphisms of the four torus, $\mathbb{T}^4$, in \cite{BV}. 
In this case, $\mathbb{T}^4=U$ and every non-trivial 
dominated splitting of the system is of the form $E\oplus F$ where both $E$ and $F$ are two-dimensional and non-hyperbolic.

We emphasize that the arguments in \cite{DN} are all two dimensional and that their translation to higher dimensions is not straightforward. 
 Most of Asaoka's paper deals with the existence of persistent tangencies and 
 the proof of non-existence of symbolic extensions is quite brief and refers to \cite{DN}.  Moreover, it considers a 
very specific dynamical configuration. 
So we consider  important to give a complete explanation of these arguments, provided in Section~\ref{s.nonexistence},  in the general case.

Theorems~A and B  raise the following natural question.

\begin{que}
If a partially hyperbolic diffeomorphism has a center bundle that splits into 1-dimensional subbundles, then does it have a principal symbolic extension?
\end{que}

\section{background}\label{s.background}

In this section we review some facts on subshifts, entropy, asymptotic $h$-expansivity, hyperbolicity, and weak forms of hyperbolicity. 

\subsection{Entropy}
If we let $\mathcal{A}=\{0,...,m-1\}$, then the {\it full $m$-shift} is the space $\Sigma_m=\mathcal{A}^{\mathbb{Z}}$ endowed with the product topology together with the map $\sigma:\Sigma_m\rightarrow \Sigma_m$ defined by $\sigma(s)=t$ where $t_i=s_{i+1}$ for all $i\in\mathbb{Z}$ and all $s=(s_i)\in\Sigma_m$.  A {\it shift space} is a closed, shift invariant subset of a full shift.

The existence and non-existence of symbolic extensions is related to the entropy structure of a dynamical system as defined by Downarowicz~\cite{Dow}.  We now review some basic definitions of topological and measure theoretical entropy.  Let $(X,d)$ be a compact metric space and $f$ be a continuous self-map of $X$.  
The $d_n$ metric on $X$ is defined as
$$d_n(x,y)=\max_{0\leq i\leq n-1}d(f^i(x), f^i(y))$$
and is equivalent to $d$ and defined for all $n\geq 0$.  
For a set $Y\subset X$ a set $A\subset Y$ is {\it $(n,\epsilon)$-spanning} if for any $y\in Y$ there exists a point $x\in A$ where $d_n(x,y)<\epsilon$.  The minimum cardinality of an $(n,\epsilon)$-spanning set of $Y$ is denoted $r_n(Y,\epsilon)$.
The {\it topological entropy} for a system $(X,f)$ is
$$h_{\mathrm{top}}(f)=\lim_{\epsilon\rightarrow 0}\, \Big(\limsup_{n\rightarrow\infty}\frac{1}{n}\log r_n(X, \epsilon) \Big).$$

Let $(X,\mathcal{B}, \mu)$ be a Lebesgue measure space with $\mu(X)=1$.  A {\it partition} of $X$ is a collection $\xi$ of measurable sets, $\xi=\{C_{\alpha}\in\mathcal{B}\, |\,\alpha\in I\}$, such that 
\begin{itemize}
\item $\mu(X/\bigcup_{\alpha\in I}C_{\alpha})=0$, 
\item $\mu(C_{\alpha})>0$ for all $\alpha\in I$,  and 
\item $\mu(C_{\alpha_1}\cap C_{\alpha_2})=0$ for $\alpha_1\neq \alpha_2$.
\end{itemize}

For partitions $\xi$ and $\nu$ the {\it joint partition of $\xi$ and $\nu$} is
$$\xi\vee\nu=\{C\cap D\, |\, C\in\xi,\, D\in\nu,\, \mu(C\cap D)>0\}.$$

Let $f$ be a measure preserving transformation of $(X,\mathcal{B}, \mu)$.  For a measurable partition $\xi$ and $n\in\mathbb{N}$ we define the {\it joint partition of $\xi$ with respect to $f$ for $n$} to be the partition

\begin{equation}\label{e.partition}
\xi^f_{n}=\xi\vee f^{-1}(\xi)\vee\cdots\vee f^{-n+1}(\xi).
\end{equation}

The {\it metric entropy of $f$ relative to the partition $\xi$} is

\begin{equation}\label{e.metric}
h_{\mu}(f, \xi)=\lim_{n\rightarrow\infty}\frac{1}{n}H_\mu (\xi^f_{n})
\end{equation}
where
$$
H_{\mu}(\xi)=-\sum_{\alpha\in I}\mu(C_{\alpha})\log\mu(C_{\alpha}).
$$
Then the {\it entropy of $f$ with respect to $\mu$} is
$$h_{\mu}(f)=\sup\{h_{\mu}(f,\xi)\, |\,\xi\textrm{ is a measurable partition with }H_\mu(\xi)<\infty\}.$$

We now review the definition of asymptotic $h$-expansivity.  
Given a subset
$Y\subset X$ we let 
$$
\bar{r}(Y,\epsilon)=\limsup_{n\rightarrow\infty}\frac{1}{n}\log r_n(Y,\epsilon) \quad \textrm{and} \quad \tilde{h}(f,Y)=\lim_{\epsilon\rightarrow 0}\bar{r}(Y,\epsilon).$$
We denote the closed ball with  center at $x$ and radius $\epsilon$ in the $d_n$ metric as $B_{\epsilon}^n(x)$.  Let 
$$\Phi_{\epsilon}(x) =\bigcap_{n=1}^{\infty}B^n_{\epsilon}(x).$$
Then 
$$h^*_f(\epsilon)=\sup_{x\in X}\tilde{h}(f,\Phi_{\epsilon}(x)).$$ 
The map $f$ is {\it asymptotically $h$-expansive} if 
$$\lim_{\epsilon\rightarrow 0}h^*_f(\epsilon)=0.$$

\subsection{Hyperbolicities}
We now review some definitions and facts of weak forms of hyperbolicity.
For a diffeomorphism $f$ of $M$  a compact $f$-invariant set $\Lambda$
has a {\it dominated splitting} if 
$$
T_\Lambda M=E_1\oplus\cdots\oplus E_k
$$
 where each $E_i$ is non-trivial and $Df$-invariant for $1\leq i\leq k$ and there exists an $m\in\mathbb{N}$ such that 
$$\|Df^n|_{E_i(x)}\|\, \|(Df^n|_{E_j(x)})^{-1}\|\leq \frac{1}{2}$$
for every $n\geq m$, $i>j$, and $x\in \Lambda$.  

The set $\Lambda$ is {\it partially hyperbolic} if it has a dominated splitting 
$$
T_\Lambda M=E_1\oplus\cdots\oplus E_k
$$ 
and there exists some $n\in\mathbb{N}$ such that $Df^n$ either uniformly contracts $E_1$ or uniformly expands $E_k$. 
A diffeomorphism $f$ is {\it partially hyperbolic} if the whole manifold $M$ is partially hyperbolic.

Finally, the set $\Lambda$ is {\it (uniformly) hyperbolic} if it has a dominated splitting $T_\Lambda M=E^s\oplus E^u$ 
where  there exists some $n\in\mathbb{N}$ such that $Df^n$ uniformly contracts $E^s$ and uniformly expands $E^u$. 

\begin{rem}[Appendix B.1 in \cite{BDV}] 
The properties of dominated splitting and partial and uniform hyperbolicity are robust in $\mathrm{Diff}(M)$.
\label{r.robustspltting}
\end{rem}

A set $\Lambda$ is {\it locally maximal} if there exists an open set $U$ containing $\Lambda$ such that 
$\Lambda=\bigcap_{n\in\mathbb{Z}}f^n(U).$  A set is {\it transitive} if there exists a point in the set whose forward orbit is dense in the whole set.  A hyperbolic set is a {\it basic set} if it is locally maximal and transitive.
One of the important properties of hyperbolic sets is structural stability. This allows one to define the continuation of a hyperbolic
sets.  See~\cite[p.  571-572]{KH}.
For a hyperbolic set $\Lambda=\Lambda_f$ of a diffeomorphism $f$ we denote the {\it continuation\/} of $\Lambda$
for $g$ close to $f$ by $\Lambda_g$.

For a partially hyperbolic diffeomorphism $f$ with a splitting $E^s\oplus E^c\oplus E^u$ we know there exist unique families $\mathcal{F}^u$ and $\mathcal{F}^s$ of injectively immersed submanifolds such that $\mathcal{F}^i(x)$ is tangent to $E^i$ for $i=s,u$, and the families are invariant under $f$.   These are called, respectively,  the unstable and stable foliations of $f$.  For the center direction it is  known, in the general case, 
that there is no foliation tangent to the center bundle \cite{Gour}.  If the center bundle is 1-dimensional, then there always exist curves tangent to the center direction through every point.  Such curves are called  {\it central curves}.  However, these curves need not be unique due to the fact that the center bundle, in general, is not Lipschitz.

\section{Partially hyperbolic diffeomorphisms with 1-dimensional center foliations}

In this section we prove Theorem~A.  We will see that the existence of 
a 1-dimensional center bundle will allow us to use 1-dimensional arguments to compute the entropy. 
Before proceeding we define some notation.  For a curve $\eta$ in $M$ its lenght is denoted by $|\eta|$.

\begin{figure}[htb]
\begin{center}
\psfrag{x}{$x$}
\psfrag{b}{$y'$}
\psfrag{v}{$W^s_{\alpha}(y)$}
\psfrag{y}{$y$}
\psfrag{u}{$z$}
\psfrag{A}{$B_{\epsilon}(x)$}
\psfrag{w}{$\gamma(x)$}
\psfrag{z}{$W^u_{\alpha}(y)$}
\psfrag{e}{$\epsilon$}
\psfrag{c}{$V^s_{\gamma(x), \alpha}(x)$}
\psfrag{d}{$V_{\gamma(x)}(x)$}
\includegraphics{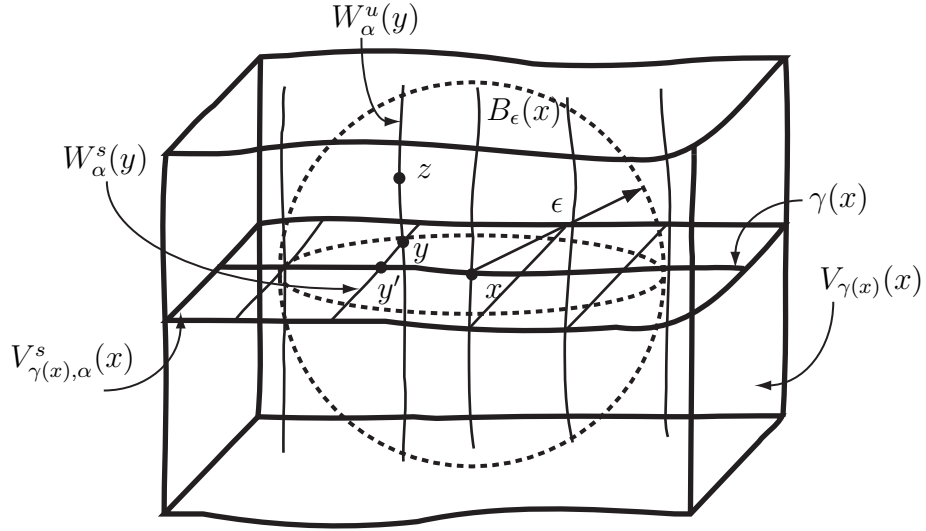}
\caption{Foliation box}\label{f.foliations}
\end{center}
\end{figure}

\subsection{Proof of Theorem~A}\label{ss.texistence}
Let $f\in\mathrm{Diff}(M)$ be partially hyperbolic with 1-dimensional center bundle.  Throughout the rest of the proof we assume that $TM=E^s\oplus E^c\oplus E^u$ with each of these is nontrivial.  The case where $E^s$ or $E^u$ is trivial is a simplification and is omitted.  Furthermore, we let $\mathrm{dim}(E^s)=s$.

The first step of the proof is to show that there exists a ``foliation chart" on a uniform scale.  More precisely,  define $\lambda$ as the maximal expansion of the derivative of $f$ and select constants $\delta_0,\delta_1, \delta_2,>0$ satisfying the following:
\begin{itemize}
 \item 
for every curve $\gamma$ tangent to $E^c$ with $|\gamma|<\lambda\, \delta_0$ and every $s$-disk $\Delta^{s}\subset \mathcal{F}^{s}$ with 
$s$-diameter less than 
$\delta_1$ one has that $\gamma \cap \Delta^{s}$ contains at most one point; and
\item
for every ($s+1$)-disks $\Upsilon$ of the form
$$\Upsilon=\bigcup_{y\in\gamma}\mathcal{F}^s_{\delta_1}(y),$$ 
where $\gamma$ is a curve tangent to $E^c$ with length bounded by $\delta_0$
and $\mathcal{F}^s_{\delta_1}(y)$ is the $\delta_1$-disk centered at $y$ contained in $\mathcal{F}^s$, and
every $u$-disk $\Delta^{u}\subset \mathcal{F}^{u}$ with $u$-diameter less than $\delta_2$ one has that $\Upsilon \cap \Delta^{u}$ 
contains at most one point.
\end{itemize}

We select small $\delta>0$ such that $\lambda \, \delta<\delta_i$ for $i=0, 1,2$.

Fix $\alpha\in(0,\delta/2)$.  For each $x\in M$ consider a curve $\gamma(x)$ tangent to $E^c$
containing $x$ with radius $\alpha$. (So the extreme points of $\gamma(x)$ are at distance $\alpha$ from $x$, where $\alpha<\delta/2$.)  This is possible since $E^c$ is
1-dimensional. 
By construction, $\mathcal{F}^{s}_{\delta_1} (z)\cap \gamma(x)$ 
consists of at most one point for all $x$ and $z$.

We let
$$
V^s_{\gamma(x),\alpha}(x)=\bigcup_{y\in \gamma(x)} \mathcal{F}^{s}_\alpha (y).
$$
For small $\tau>0$ we define
$$
V^{s,u}_{\gamma(x),\alpha,\tau}(x)=\bigcup_{z\in V^s_{\gamma(x),\alpha(x) }} \mathcal{F}^{u}_\tau (z).
$$
For notational simplicity, we fix small $\alpha>0$ and write
$$
V_{\gamma(x)}(x)=V^{s,u}_{\gamma(x),\alpha,\alpha}(x).
$$

For every $\epsilon>0$ sufficiently small and  $x\in M$
we have
$$
B_{\epsilon}(x)\subset V_{\gamma(x)} (x).
$$
Note that this inclusion holds for all $\gamma(x)$ as above.
See Figure~\ref{f.foliations}.  

\begin{rem}[Foliation chart]
We selected $\delta_0$ in such a way that for all central curve $\eta$ of size $\lambda\, \delta_0$
the cube
 $$
V_{\eta}=V^{s,u}_{\eta,\alpha,\alpha}
$$
provides a foliation chart.  (See for instance~\cite[p. 19]{CC} for a definition of a foliation chart.) 
\label{r.foliationchart}
\end{rem}

\begin{lem}
For every $x\in M$ and $\epsilon$ small enough it holds that
$$
\Phi_\epsilon(x)= \bigcap_{n=1}^\infty B_\epsilon^n(x)\subset V_{\gamma(x),\alpha}^s (x).
$$
\label{l.central}
\end{lem}

\begin{demo}
We consider the forward orbit of $x$ and define $x_n= f^n(x)$.  For each $n$ we define a central  curve 
$\gamma_n =\gamma (f^n(x))=\gamma(x_n)$ such that $\gamma_n$ is ``compatible'' with
$f(\gamma_{n-1})$: this means we consider the components $\gamma^\pm_n (x_n)$ of $\gamma_n(x_n)\setminus \{x_n\}$ and have
either 
\begin{itemize}
\item $f(\gamma_n^\pm) \subset \gamma_{n+1}^\pm$ or 
\item $\gamma_{n+1}^\pm \subset f(\gamma_{n}^\pm)$ or 
\item $f(\gamma_n^\pm) \cap \gamma_{n+1}^\pm=\emptyset$.
\end{itemize}
  
The ``cubes'' $V_{\gamma_n}(x_n)$ containing $B_\epsilon(x_n)$ are defined as above.


Let $z\in \Phi_{\epsilon}(x)$. 
By definition and construction, $f^n(z)=z_n \in V_{\gamma_n}(x_n)$ for all $n\ge 0$.
We let $y$ be the unique point in $V^s_{\gamma(x), \alpha}(x)$ such that $z\in \mathcal{F}^u_{\alpha}(y)$ and $y'$ be the unique point in $\gamma(x)$ such that $y\in \mathcal{F}^s_{\alpha}(y')$.  See Figure~\ref{f.foliations}.

We define $y_n$, $y_n^\prime$, and $z_n$ inductively. 
By construction, $y_n^\prime\in \gamma_n$, otherwise, by Remark~\ref{r.foliationchart}, $z_n\notin B_\epsilon (x_n)$.

If $z_n\neq y_n$, then the uniform expansion along the unstable direction implies that there is a first $n$ such that $z_n$ and $y_n$ can not be
in the same foliation chart of radius $\alpha$.  
To be more precise, by hypothesis, we know that $y_n$ is in $V^s_{\gamma_n}(x_n)$ for all $n$.
If $z_n\ne y_n$ for some $n$ the uniform expansion
along the unstable direction implies that there is a first $m$ 
such that
$z_m\not\in V_{\gamma_m}(x_m)$.  Hence, $z_m\notin B_{\epsilon}(x_m)$, a contradiction. 
This ends the proof of the lemma.
\end{demo}

Define
$$
\Gamma^c(x)= \bigcap_{n\ge 0} f^{-n} (\gamma_n).
$$

Given a central curve $\eta$ of size at most $\delta_0$ and a point $y\in V^s_{\eta, \alpha}$ 
we let $y^\prime$ be the unique point in $\eta$ such that
$y\in \mathcal{F}^{s}_\alpha (y^\prime)$.

Note that the proof of Lemma~\ref{l.central} gives the following:
\begin{cla}
Let $y\in \Phi_\epsilon(x)$. Then $y^\prime \in \Gamma^c(x)$.
\end{cla}

This claim implies that
$$
|f^n (\Gamma^c(x))|< 2\, \delta_0.
$$
A folklore fact
implies that the growth of  an $(n,\delta)$-spanning set in $\Gamma^c(x)$ is
subexponential for all $\delta>0$, see  for instance \cite[Lemma 3.2]{BFSV}.

The next lemma will show that $f$ is asymptotically $h$-expansive, completing the proof of Theorem~A.

\begin{lem}\label{l.hexpansive}
For every $x\in M$ and all $\delta>0$ small enough, 
$$
\tilde{r}(\Phi_{\epsilon}(x), \delta)=\limsup_{n\rightarrow\infty}\frac{1}{n}r_n(\Phi_{\epsilon}(x),\delta)=0.$$  
\end{lem}

\begin{demo}
 Since $f$ is uniformly contracting in the stable direction we know that any exponential growth of the spanning sets for 
$\Phi_{\epsilon}(x)$ occurs along $\Gamma^c (x)$. This means that we can focus on the growth of 
$(n,\delta)$-spanning set in $\Gamma^c(x)$.  More
precisely: 
\begin{fac}
Let $Y$ be an $(n,\delta/2)$-spanning set in $\Gamma^c(x)$. Then $Y$ is a $(n,\delta)$-spanning set in  
$V^s_{\Gamma^c(x), \alpha}$ for all $n$ sufficiently large.
\end{fac}
By the comments above, this growth is subexponential and the proof of the theorem is complete.
\end{demo}

\subsection{Proof of Corollary~\ref{c.existence}}\label{ss.cexistence}
It is enough to argue exactly as in Subsection~\ref{ss.texistence} considering a small neighborhood of the
robustly transitive set where the bundles can be extended (recall Remark~\ref{r.robustspltting}) and consider balls in this neighborhood.

\section{Diffeomorphisms with no symbolic extensions}\label{s.nonexistence}

This section consists of three parts. In the first one, 
in the context of non-hyperbolic robustly transitive sets (and diffeomorphisms),
we explain how 
the existence of a center indecomposable bundle of dimension at least 2 
generates persistent homoclinic tangencies, see Proposition~\ref{p.bundles}.
To get the non-existence of symbolic extensions we will use the 2-dimensional arguments in \cite{DN}. 
For that we need to consider
special saddles in the robustly transitive set (saddles with real multipliers) and to see that the tangencies
occur in a local normally hyperbolic surface.
In the second one, we explain how the non-existence of symbolic extensions follows from this property. Finally, 
in the last section we prove
Theorem~B.

\subsection{Persistence of homoclinic tangencies} 
In this section we consider diffeomorphisms in 
  $\mathcal{T}^{nh}_2(U)$, that is, the subset of $\mathrm{Diff}(M)$ of 
diffeomorphisms $f$ such that $\Lambda_f(U)=\bigcap_{n\in \mathbb{Z}} f^n(\bar U)$ is a robustly transitive set whose
finest dominated splitting 
has some (indecomposable)  center bundle with dimension at least 2. We begin by reviewing the 
constructions in \cite{BoDiAENS,BDP,BDPR,BoDipreprint}.

Define $\mathcal{N}(U)$ as the set of diffeomorphisms $f\in \mathrm{Diff}(M)$ such that $\Lambda_f(U)$ is robustly transitive and
contains periodic points of different indices. Here the {\it index\/} of a periodic point is the dimension of its stable bundle.
We let $s^-(f)$ and $s^+(f)$ be the minimum and maximum of the indices of periodic points
of $\Lambda_f(U)$.

The set $\mathcal{N}(U)$ is open and dense in the set $\mathcal{T}^{nh}(U)$ of diffeomorphisms $f$ such that $\Lambda_f(U)$ is non-hyperbolic and
robustly transitive, see \cite{BDPR}.

Given
a diffeomorphism 
$f\in \mathcal{N} (U)$, 
we consider the finest dominated splitting 
$E_1\oplus E_2 \oplus \cdots \oplus E_{k(f)}$ of $\Lambda_f(U)$, see \cite[Theorem 4]{BDP}.
This means that every subbundle $E_i$ of the splitting  is indecomposable.
We let $\alpha=\alpha(f)$  be the first $j\in \{1,\dots,k(f)\}$ such that 
$E_j$ is not uniformly contracting.
Similarly, we let $\beta$ be the last $j\in \{1,\dots,k(f)\}$  such that $E_j$ is not uniformly expanding.
Then $E_1\oplus \cdots \oplus  E_{\alpha-1}$ is a uniformly contracting bundle and 
$E_{\beta+1}\oplus \cdots \oplus   E_{k(f)}$ is uniformly expanding.

\begin{rem}[Corollary C in \cite{BDPR}] \label{r.bdpr}
Define $d_i=\dim (E_i)$. There is an open and dense subset $\mathcal{V}(U)$ 
of $\mathcal{N}(U)$ such that following maps are locally constant: 
the number $k(f)$ of bundles of the finest dominated splitting,
the dimensions $d_i$ of these bundles, and the numbers $\alpha(f)$ and $\beta(f)$. 

Furthermore, the functions $s^-(f)=s^-$ and $s^+(f)=s^+$ are locally constant in $\mathcal{V}(U)$ and for every $j\in [s^-,s^+]$ 
there is a saddle of index $j$ in $\Lambda_f(U)$. Finally, $s^-=d_1+\cdots + d_{\alpha-1}$ and $s^+=d_1+\cdots + d_{\beta}$.
\end{rem}

We are ready to state a key result about  persistence of tangencies.

\begin{prop}\label{p.bundles}
Let $f\in\mathcal{V}(U)\cap\mathcal{T}^{nh}_2(U)$ and let $\mathbb{E}_1 \oplus\cdots \oplus\mathbb{E}_k$ be the finest dominated splitting of $\Lambda_f(U)$. Consider 
$l\in[\alpha(f),\beta(f)]$ such that $d_l\geq 2$.  

Suppose that $p_f\in \Lambda_f(U)$ is a saddle 
of index $d_1 +\cdots +d_{l-1}+1$ such that its continuation is defined for every $g$ in a neighborhood $\mathcal{U}_f\subset \mathcal{V}(U)\cap\mathcal{T}^{nh}_2(U)$ of $f$.  Then there is a dense subset $\mathcal{T}$ 
of $\mathcal{U}_f$ consisting of diffeomorphisms with homoclinic tangencies associated to $p_g$.
\end{prop} 

Note that if $f\in\mathcal{V}(U)\cap\mathcal{T}^{nh}_2(U)$ then there exists some $l\in[\alpha,\beta]$ such that $d_l\geq 2$ and therefore there is a saddle $p_f$ as in  Proposition~\ref{p.bundles}. 

We observe that the conclusion of the proposition also holds
if the saddle  $p_f\in \Lambda_f(U)$ 
has index $\ell=d_1 +\cdots +d_{l-1}+j$, $j\ge 1$, with  $\ell<d_1 +\cdots +d_l$.

In what follows, we fix a component $\mathcal{C}(U)$ of  $\mathcal{V}(U)\cap\mathcal{T}^{nh}_2(U)$ where the maps
in Remark~\ref{r.bdpr} are constants and let 
$$
r=d_1 +\cdots +d_{l-1}+1.
$$

We need to introduce some notation. 
Given a periodic point $q$ we denote by $\tau (q)$ its period and 
by 
$$\lambda_1(q), \dots ,\lambda_n(q)$$ 
the eigenvalues of $Df^{\tau (q)}(q)$ counted with multiplicity and ordered in non-decreasing modulus.
We say that $\lambda_i(q)$ is the {\it $i$-th eigenvalue\/} of $q$.

We define $\mathrm{Per}_{\mathbb{R}}(f,U)$ as the set of saddles $q\in \Lambda_f(U)$ such that
all eigenvalues of $q$ are real and different in modulus.
We let  $\mathrm{Per}_\mathbb{R}(f,U)^{j}$ be the subset of $\mathrm{Per}_{\mathbb{R}}(f,U)$ of saddles of index $j$.
We also consider the subset
$\mathrm{Per}(f,U)^{j}_n$ of $\mathrm{Per}_{\mathbb{R}}(f,U)^{j}$ of points with period $n$. 

We define $\mathrm{Per}(f,U)_{\mathbb{C}}^{r,r+1}$ as the set of saddles $q\in \Lambda_f(U)$ such that
$$
|\lambda_{r-1}(q)|<|\lambda_{r}(q)|=|\lambda_{r+1}(q)|<1<|\lambda_{r+2}(q)|
$$ 
and $\lambda_{r}(q)$ and $\lambda_{r+1}(q)$ are non-real. Note that $q$ has index $r+1$.

Fixed $r$ as above
consider the following subsets of $\mathcal{C}(U)$:
\begin{itemize}
 \item 
$\mathcal{D}(U) : = \{g\in \mathcal{C}(U) \colon $ 
$\mathrm{Per}_{\mathbb{R}}(g,U)^{r}$  and $\mathrm{Per}_{\mathbb{R}}(g,U)^{r+1}$ are both dense in $\Lambda_g(U)\}$.
\item
$\mathcal{F}(U) : = \{g\in \mathcal{C}(U) \colon $ such that 
$\mathrm{Per}_{\mathbb{C}}(g,U)^{r,r+1}$ is open and dense in $\mathcal{C}(U)\}$.
\item
$\mathcal{O}(U): = \{ g\in \mathcal{C}(U)\colon$ the sets $\mathrm{Per}_{\mathbb{R}}(g,U)^{r}$,
$\mathrm{Per}_{\mathbb{R}}(g,U)^{r+1}$ and 
$\mathrm{Per}_{\mathbb{C}}(g,U)^{r,r+1}$  are all non-empty$\}$.
\end{itemize}

Next remark is just a dynamical reformulation of some results about cocyles
(Lemmas 4.16 and 5.4) in \cite{BDP}.

\begin{rem}\label{r.index}
Let $\mathcal{C}(U)$ and $r$ be as above. 
\begin{itemize}
 \item
$\mathcal{D}(U)$ is residual in $\mathcal{C}(U)$.
\item
$\mathcal{F}(U)$ and $\mathcal{O}(U)$ are open and dense in $\mathcal{C}(U)$.
\end{itemize}
\end{rem}

An immediate and standard consequence of Hayashi's Connecting Lemma \cite{Hayashi} is the next remark.

\begin{rem}\label{r.connecting} There is a residual subset of $\mathcal{C}(U)$ such 
that for every $g$
and any pair of saddles $q_g$ and $q_g^\prime$ in $\Lambda_g(U)$ of indices $r+1$ and $r$
one has that
$W^s(q_g)$ and $W^u(q^\prime_g)$ have some transverse intersection. Moreover, there are $h$ arbitrarily close to $g$ such that
$W^u(q_h)\cap W^s(q^\prime_h)\ne \emptyset$.
\end{rem}

\begin{figure}[htb]
\begin{center}
\psfrag{a}{$p_g$}
\psfrag{b}{$q_g$}
\psfrag{c}{$q_{\bar{g}}$}
\psfrag{d}{$W^u(p_{\bar{g}})$}
\psfrag{e}{$W^s(p_{\bar{g}})$}
\includegraphics{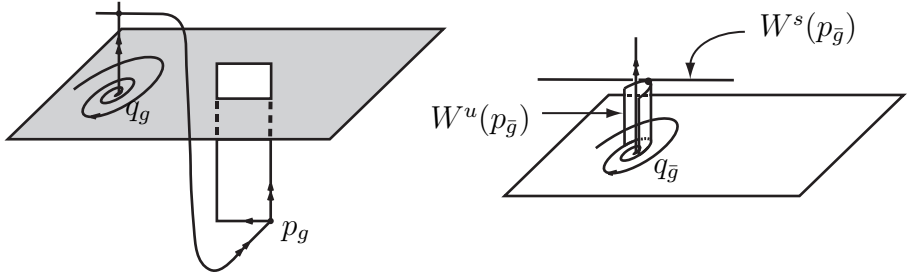}
\caption{Unfolding of a cycle for $p_g$ to create a tangency}\label{f.tangency}
\end{center}
\end{figure}

\subsubsection{Sketch of the proof of Proposition~\ref{p.bundles}}
If we omit the condition of the persistence of the tangencies, this proposition is just Theorem F in \cite{BDPR}.
We now explain how the
persistence is obtained. This follows from standard arguments in the $C^1$-topology
that can
be found in several papers. We now recall the main steps.

Fix a saddle $p_g$ of index $r$ as in the proposition.
By Remark~\ref{r.index}, after a perturbation  
we can assume that $\Lambda_g(U)$ contains a saddle 
$q_g\in \mathrm{Per}_{\mathbb{C}}(g,U)^{r,r+1}$ 
 of index
$r+1$.
By Remark~\ref{r.connecting}, we can assume that $W^u(p_g)$ and $W^s(q_g)$ has some transverse
intersection. Therefore, we can assume that $W^u(p_g)$ spiralizes around $W^u(q_g)$. By Remark~\ref{r.connecting}, after a perturbation,
we can get an intersection between  $W^s(p_g)$ and $W^u(q_g)$ (obtaining a 
cycle associated to $p_g$ and $q_g$). Unfolding this cycle, one gets a homoclinic tangency assoicated to $p_g$.
These arguments are depicted in Figure~\ref{f.tangency}. See
for instance
Section 2.2 \cite{BoDiAENS} for details.

\subsection{Downarowicz-Newhouse construction}

To begin this section we review the condition given in~\cite{DN} that guarantees the non-existence of symbolic extensions, see Theorem~\ref{t.dn}.  Using Proposition~\ref{p.bundles} we will then verify that there is a residual set $\mathcal{R}$ in 
$\mathcal{T}^{nh}_{2}(U)$ such that each $f\in\mathcal{R}$ satisfies this  condition, and therefore every $f\in \mathcal{R}$ has no
symbolic extension.

As we will examine invariant measures for a system $(X, f)$, we review some basic facts of invariant measures.  Let  
$\mathcal{M}(f, X)$ be the set of invariant measures  and $\mathcal{M}_e(f, X)$ be the set of ergodic invariant measures for $(X,f)$.  If $f$ is  a homeomorphism of a compact metric space then  $\mathcal{M}(f, X)$ is nonempty and satisfies the following (see~\cite[Section 4.6]{BS}):
\begin{itemize}
\item $\mathcal{M}(f, X)$ is convex and compact for the weak$^*$-topology, and
\item $\mathcal{M}_e(f, X)$ is precisely the extreme points of $\mathcal{M}(f, X)$.
\end{itemize} 
Let $\rho(\cdot, \cdot)$ denote the metric on $\mathcal{M}(f, X)$ giving the weak$^*$-topology.

A sequence of partitions $\{\alpha_k\}$ is {\it essential} if
\begin{itemize}
\item $\mathrm{diam}(\alpha_k)\rightarrow 0$ as $k\rightarrow 0$, and
\item $\mu(\partial\alpha_k)=0$ for any $\mu \in \mathcal{M}(f, X)$, where $\partial \alpha_k$ is the union of the boundaries of the elements in the partition $\alpha_k$.
\end{itemize}

A {\it simplicial sequence of partitions} is a sequence $\{\alpha_k\}$ of nested partitions (each $\alpha_k$ is a refinement of $\alpha_{k-1}$) whose diameters go to zero and such that for each $\alpha_k$ the partition is given by a smooth triangulation of $M$.

\begin{prop}[Proposition 4.1 in \cite{DN}]
\label{p.simplicial}
Let $\{\alpha_k\}$ be a simplicial sequence of partitions of $M$.  Then there is a residual set $\mathcal{S}\subset \mathrm{Diff}(M)$ such that if $f\in\mathcal{S}$, then $\{\alpha_k\}$ is an essential partition for $f$.
\end{prop}

The next result is the condition used  to prove that there is no symbolic extension.

\begin{thm}[Proposition 4.4 in \cite{DN}]
\label{t.dn}
Let $\{\alpha_k\}$ be an essential sequence of partitions on $M$.  Suppose there exists an $\epsilon>0$ and a compact set $\mathcal{E}\subset\mathcal{M}(f, M)$ such that for every $\mu\in \mathcal{E}$ and every $k>0$, 
$$\limsup_{\nu\rightarrow\mu, \nu\in\mathcal{E}}\, \big(h_{\nu}(f)-h_{\nu}(\alpha_k, f)\big)>\epsilon.$$
Then  a symbolic extension for $f$  does not exist.
\end{thm} 

The construction involves measures supported on hyperbolic basic sets. By \cite{Sig} these measures can be obtained as limits of
periodic measures defined as follows.
Let $f\in\mathrm{Diff}(M)$ and $p$ be a hyperbolic periodic point for $f$ with period $\tau(p)$.  Then the {\it periodic measure for $p$} is
$$\mu_p=\frac{1}{\tau(p)}\sum_{i=0}^{\tau(p)-1}\delta_{f^i(p)},$$ where $\delta_{f^i(p)}$ is the point mass at $f^i(p)$.  

Using Proposition~\ref{p.bundles} we will construct $\mathcal{E}$ as the closure of a certain set of periodic measures.  
 
Note that it is enough to prove Theorem~B for the set
$\mathcal{V}(U) \cap \mathcal{U}_f$, see
Remark~\ref{r.bdpr} and Proposition~\ref{p.bundles}. Therefore the number 
$$r=d_1+\cdots +d_{l-1}+1$$  
in Propostion~\ref{p.bundles} remains fixed from now on.

Let $\lambda_j(q)$ be the $j$-th eigenvalue of a saddle  $q$. We let $\chi_j(q)$ be the $j$-th 
Lyapunov exponent of $q$, i.e., 
$$\chi_j(q)= \frac{1}{\tau (q)}\, \log (|\lambda_j(q)|),$$ and define 

\begin{equation}\label{e.inf}
\chi_{j} (f,U)=\inf \{\chi_{j} (q)\, :\, q\in \mathrm{Per}_{\mathbb{R}}(f,U)^{j}\}.
\end{equation}
For an ergodic measure $\mu$  we let  $\chi_j(\mu)$ be its $j$-th Lyapunov exponent.

A hyperbolic set $\Lambda$ is {\it subordinate} to a partition $\xi$ if there exist compact sets $\Lambda_1,...,\Lambda_j$ for some $j\geq 1$ such that 
\begin{itemize}
 \item $\Lambda=\bigcup_{i=1}^j\Lambda_i$,
\item $f(\Lambda_i)=\Lambda_{i+1}$, for $1\leq i\leq j-1$, and $f(\Lambda_j)=\Lambda_1$, and
\item $\Lambda_i$ is contained in a single element of the partition $\xi$ for all $1\leq i\leq j$.
\end{itemize}

\begin{defn} A diffeomorphism $f$ with an essential sequence of simplicial partitions $\{\alpha_k\}$ satisfies property $\mathcal{S}^j_n$ if for 
each $q\in \mathrm{Per}_{\mathbb{R}}(f,U)^{j}_n$ the following hold: 
\begin{enumerate}
\item[(a)] there exists a 0-dimensional hyperbolic basic set $\Lambda(q,n)$ for $f$ such that 
$$
\Lambda(q,n)\cap \partial \alpha_n=\emptyset\textrm{ and } \Lambda(q,n)\subset \Lambda (f,U),$$
\item[(b)] the set $\Lambda(q,n)$ is subordinate to $\alpha_n$, 
\item[(c)] there is a $\mu\in\mathcal{M}_e(f, \Lambda(q,n))$ such that 
$$|h_{\mu}(f)-\chi_j(q)|<\frac{1}{n}\chi_j(q), $$
and
\item[(d)] for every $\mu\in \mathcal{M}_e(f, \Lambda(q,n))$ we have
$$\rho(\mu, \mu_q)<\frac{1}{n}\textrm{ and }|\chi_j(\mu)-\chi_j(q)|<\frac{1}{n}\chi_j(q).$$
\end{enumerate}
\label{d.essential}
\end{defn}

Let $f\in \mathcal{O}(U)$ (see Remark~\ref{r.index})
and let $\tau_r(f,U)$ be the minimum period of points in 
$\mathrm{Per}_{\mathbb{R}}(f,U)^{r}$. 
We let $\mathcal{R}^r_{m}$ be the set of diffeomorphisms in $\mathcal{O}(U)$ 
with $\tau_r(f,U)=m$.

For $m\leq n$ let $\mathcal{D}^r_{m,n}$ be the set of $\mathcal{R}^r_{m}$ of diffeomorphisms satisfying property $\mathcal{S}^r_n$.

\begin{lem}\label{l.opendense}
For every $m$ and sufficiently large $n$ the set $\mathcal{D}^{r}_{m,n}$ is open and dense in $\mathcal{R}^{r}_{m}$.
\end{lem}

\begin{demo}
The proof of the above lemma follows from the proof of Lemma 5.1 in~\cite{DN}. 
Since the sets $\Lambda(q,n)$ are hyperbolic the sets $\mathcal{D}^r_{m,n}$ are open  (see \cite[page 471]{DN}). Thus the 
crucial point is density and for that we use Proposition~\ref{p.bundles}. 

\begin{figure}[htb]
\begin{center}
\psfrag{q}{$q$}
\psfrag{p}{$q$}
\psfrag{c}{$c$}
\psfrag{d}{$\pi$}
\psfrag{e}{$\pi$}
\includegraphics{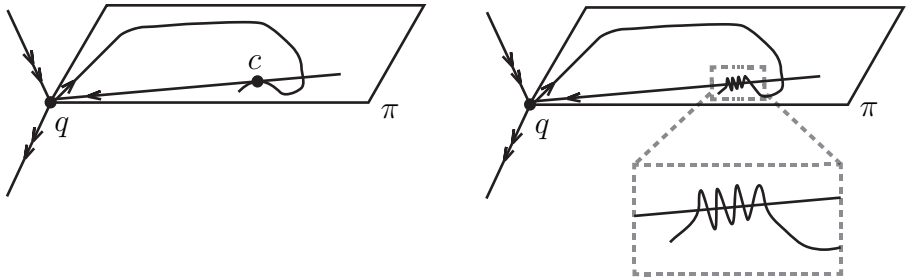}
\caption{The creation of $\Lambda(q,n)$ in the planar region $\pi$}\label{f.three}
\end{center}
\end{figure}

Fix $n$ and $m$, with $m\le n$,
and consider a diffeomorphism $g\in\mathcal{R}_m^{r}$ and a saddle $q\in \mathrm{Per}_{\mathbb{R}}(g,U)^{r}_n$.
By Proposition~\ref{p.bundles}, after arbitrarily small perturbation we can assume that $g$ has a
 homoclinic tangency associated to $q$. Since the  $r$-th and $(r+1)$-th multipliers of $q$ satisfy 
$$
|\lambda_{r-1}(q)|<|\lambda_{r}(q)|<1<|\lambda_{r+1}(q)|<|\lambda_{r+2}(q)|,
$$ 
after a new perturbation,  we
can assume that this tangency occurs in a ``normally hyperbolic surface'' $\pi$: so there are 
curves $\gamma^s$ and $\gamma^u$ contained in
$W^s(q)\cap \pi$ and $W^u(q)\cap \pi$ having a quadratic homoclinic tangency.
See for instance \cite[Affirmation]{BoDiAENS} and Figure~\ref{f.three}. Therefore, from now on the arguments
are ``2-dimensional'' and we can follow \cite{DN}. 

More precisely, let $c$ be the tangency point. After a perturbation, for each large $n$, we get a hyperbolic set $\Lambda(q,n)$ of $f$ 
subordinate  to the partition $\alpha_n$ and contained in a small neighborhood of the orbit of $c$ such that
$h_{\mathrm{top}}(f,\Lambda(q,n))$ is close to $\chi_{r+1}(q)$, see \cite[equation (43), page 484]{DN}. 
By construction, the set $\Lambda(q,n)$ is contained in $\Lambda_f(U)$.
This completes the sketch of the proof of Lemma~\ref{l.opendense}.
\end{demo}

\begin{rem}\label{r.rr+1}
By construction, each periodic point $\bar q\in \Lambda(q,n)$ satisfies Proposition~\ref{p.bundles} and belong to
$\mathrm{Per}_\mathbb{R}(f,U)^{r}$. Therefore, we can apply the previous construction to these saddles.
\end{rem}

\subsection{Proof of Theorem~B}  Note that it is enough to prove a local version of this theorem for diffeomorphisms 
$g\in \mathcal{U}_f\cap \mathcal{V}_f \cap \mathcal{S}$, where $\mathcal{S}$ is the residual set in Proposition~\ref{p.simplicial}.
Recall that for these diffeomorphisms the functions $\alpha$, $\beta$, $s^-$, and $s^+$ are constant.

The proof is similar to the proof of Theorem 1.3 in~\cite[pages 471-2]{DN}. 
  Recalling equatiom \eqref{e.inf} define
\begin{equation}\label{e.rho}
\rho_0=\dfrac{\chi_{r+1}(g,U)}{2} 
\end{equation}
and let
$$
\mathcal{E}_1(g)=\{\mu_q\, :\, q\in \mathrm{Per}_{\mathbb{R}}(g,U)^r\textrm{ and }\chi_{r+1}(q)>\rho_0\},$$
and set $\mathcal{E}(g)$ as the closure of $\mathcal{E}_1(g)$.

It is sufficient to verify the conditions of Theorem~\ref{t.dn} for $\mu_q\in \mathcal{E}_1(g)$. More precisely, from Lemma~\ref{l.opendense} 
we know 
that for each $n$ sufficiently large there exists a periodic hyperbolic basic set $\Lambda(q,n)$ that is subordinate to $\alpha_n$ and all ergodic measures of $\Lambda(q,n)$ are $\frac{1}{n}$ close to $\mu_q$. Moreover,  from (d) in Definition~\ref{d.essential},
$$\chi_{r+1}(\mu)>\frac{n-1}{n}\, \chi_{r+1}(q)$$
and  from (c) in Definition~\ref{d.essential} for every $\nu_n\in\mathcal{M}_e(g, \Lambda(q,n))$ it holds that 
\begin{equation}\label{e.h}
 h_{\nu_n}(g)>\frac{n-1}{n}\, \chi_{r+1}(q).
\end{equation}

\begin{lem} \label{l.ergodic}
 $\mathcal{M}_e(g, \Lambda(q,n))\subset \mathcal{E}(g). $
\end{lem}

\begin{demo}
 Just recall that every ergodic measure supported on a hyperbolic basic set is the weak$^*$-limit of periodic measures for periodic points in the set, see \cite{Sig}. Moreover, 
by Remark~\ref{r.rr+1},
the periodic points in $\Lambda(q,n)$ belong to $\mathrm{Per}_\mathbb{R}(g,U)^r$.
\end{demo}

\begin{lem} \label{l.entropy}
Let $\nu_n\in\mathcal{M}_e(g, \Lambda(q,n))$. Then $h_{\nu_n}(\alpha_n, g)=0$. 
\end{lem}

\begin{demo}
Recall that (see equations \eqref{e.partition} and \eqref{e.metric})
$$
h_{\nu_n} (g,\alpha_n)=\lim_{j\to \infty} \, H_{\nu_n} ( (\alpha_n)_j^g)=
\lim_{j\to \infty} \, \frac{-1}{j} \, \sum_{A\in (\alpha_n)_j^g}  \nu_n(A)\, \log \nu_n(A).
$$
Since $\Lambda(q,n)$ is subordinate to $\alpha_n$ the sum is constant for each $j$, see~\cite[464]{DN}. Therefore this limit goes to
$0$ as $j\to \infty$. 
\end{demo}

Fix $k_0\in\mathbb{N}$ and 
a measure $\mu\in \mathcal{E}(g)$.  By Theorem~\ref{t.dn} it is enough  to prove that  there exists a sequence $\nu_n\rightarrow \mu$ such that 
\begin{equation}\label{e.end}
\limsup_{n\to \infty}\, \big(h_{\nu_n}(g)-h_{\nu_n}(\alpha_{k_0}, g)\big)>\rho_0.
\end{equation}
Fix a sequence $\eta_n\in \mathcal{E}(g)$ converging to $\mu$.

By definition of $\mathcal{E}(g)$, we know that for each $\eta_n$ there is is a periodic measure $\mu_q\in \mathcal{E}_1(g)$ arbitrarily close $\eta_n$. 
Then by (d) in Definition~\ref{d.essential} and Lemma~\ref{l.ergodic}
there exists
$\nu_n\in \mathcal{M}_e(g, \Lambda(q,n))\subset \mathcal{E}(g)$ which is close to $\mu_q$.

Take large  $n \ge k_{0}$, the $\alpha_n$ is a refinement of $\alpha_{k_0}$,  we have from Lemma~\ref{l.entropy} that $h_{\nu_n}(\alpha_{k_0}, g)=0$.  
Then, from \eqref{e.h} and the choice of $\rho_0$ in equation~\eqref{e.rho},  for $n$ sufficiently large we have 
$$
h_{\nu_n}(g)-h_{\nu_n}(\alpha_{k_0}, g)>\frac{n-1}{n}\chi_{r+1}(q)>\frac{\chi_{r+1}(g,U)}{2}=\rho_0.
$$
Since $\nu_n$ converges to $\mu$ we get \eqref{e.end}, ending the proof of
the theorem.
$\Box$

\bibliographystyle{plain}
\bibliography{symbolic}

\end{document}